%\hfuzz8.6pt
%\emergencystretch4mm
\input AHTOH-E.STY
%%%%%%%%%%%%%%%%%%%%%%%%%%%%%%%
\UDC{\let\,\thinspace\ 512.544.44\,+\,512.543.22\,+\,512.544.33}
\MSC{\kern-3pt 20e36,20f12,20f16,20f18,20f19,20e10}
\title
{A short proof of the Khukhro--Makarenko theorem\\
on large characteristic subgroups with laws}

\author{Anton A. Klyachko\quad{and}\quad Yulia B. Melnikova}

\address{\myAddress\\ yuliamel@mail.ru}
\grants{\RFBR 08-01-00573}

\medskip
\leftline{\it
Dedicated to A. L. Shmelkin on the occasion of his seventieth birthday
}
\bigskip

\abstract{%
We give a short proof and some strengthening of the
Khukhro--Makarenko theorem that each group virtually satisfying an
outer commutator identity contains a finite-index characteristic subgroup
satisfying this identity. An estimate for the index
of this characteristic subgroup is obtained.

{\it Key words}:
characteristic subgroups, outer commutator identities.
}
%%%%%%%%%%%%%%%%%%%%%%%%%%%%%%%%%%%%%%%%%%%%%%%%%%%%%%%%%%%%%%%%%%%%%%%%%%
%%%%%%%%%%%%%%%%%%%%%%%%%%%%%%%%%%%%%%%%%%%%%%%%%%%%%%%%%%%%%%%%%%%%%%%%%%

\noindent
Suppose that $G$ is a group and $H$ is its subgroup of finite index.
Textbooks in group theory (see, e.g., [KaM82]) contain
some simple facts allowing us to find
a finite-index subgroup in $G$ which is similar to, but better than $H$.
In particular,
\-
$H$ contains a normal finite-index subgroup of $G$
(whose index divides $|G{\,:\,}H|!$);
\-
if $G$ is finitely generated, then $H$ contains a finite-index subgroup
which is fully characteristic (and even verbal) in $G$;
\-
if $H$ is abelian, then $G$ has a characteristic abelian
subgroup of finite index.

\noindent
Recently, the last statement was substantially generalised.

\proclaim{Khukhro--Makarenko Theorem \rm([KhM07], see also [MKh07])}.
If a group $G$ has a finite-index subgroup satisfying an outer
{\rm(\emph{multilinear}, in authors' original terminology)} commutator
identity, then $G$ also has a finite-index characteristic subgroup
satisfying this identity.

Examples of outer commutator identities are $n$-step nilpotency or
solvability. The general definition is as follows. Let $F(x_1,x_2,\dots)$
be a free group of countable rank. \emph{An outer commutator of weight 1}
is a generator $x_i$.  \emph{An outer commutator of weight $t>1$} is a
word of the form
$w(x_1,\dots,x_t)=[u(x_1,\dots,x_r),v(x_{r+1},\dots,x_t)]$, where $u$ and
$v$ are outer commutators of weights $r$ and $t-r$, respectively.
Informally, an outer commutator of weight $t$ is an expression
$[x_1,x_2,\dots,x_t]$ with some arrangement of brackets.
\emph{An outer commutator identity of weight $t$} is an identity
$w(x_1,\dots,x_t)=1$ whose left-hand side is an outer commutator of weight
$t$.

The following proof of the Khukhro--Makarenko theorem is significantly
simpler and shorter than the original argument.

\bigskip
\noindent
Suppose that $H_1,\dots, H_t$ are normal subgroups of a group
$G$ and $w(x_1,\dots,x_t)$ is an outer commutator. Then
\item{1)}
the subgroup
$w(H_1,\dots,H_t)\:=\gp{w(h_1,\dots,h_t)\;;\;h_i\in H_i}$ is normal
in $G$;
\item{2)}
$w(G,\dots,G)=1$ if and only if
$G$ satisfies the identity $w(x_1,\dots,x_t)=1$;
\item{3)}
$w(H_1,\dots,H_t)=[u(H_1,\dots,H_r),v(H_{r+1},\dots,H_t)]$ if
$w(x_1,\dots,x_t)=[u(x_1,\dots,x_r),v(x_{r+1},\dots,x_t)]$;
 \item{4)}
$w(H_1,\dots,H_{i-1},\prod\limits_{N\in \cal N}\!\!N,H_{i+1},\dots,H_t)=
\prod\limits_{N\in \cal N}\!\!w(H_1,\dots,H_{i-1},N,H_{i+1},\dots,H_t)$
\newline
for any family $\cal N$ of normal subgroups of $G$.

\medskip
\enditem
These property are almost obvious and can be easily verified by
induction.

%%%%%%%%%%%%%%%%%%%%%%%%%%%%%%%%%%%%%%%%%%%%%%%%%%%%%%%%%%%%%%%%%%%%%%%%%%
%\s 1.

\Lemma.
If $w(x_1,\dots,x_t)$ is an outer commutator,
$m$ is a positive integer, $G$ is a group, and $\cal N$ is a family of
its normal subgroups such that
$$
w(\underbrace{N,N,\dots,N}_{m\rm\;times},G,G,\dots,G)=1
\quad\hbox{for all $N\in\cal N$},
$$
then
$$
w(\underbrace{\^N,\^N,\dots,\^N}_{m-1\rm\;times},\^G,\^G,\dots,\^G)=1,
\quad\hbox{where } \^N=\bigcap_{N\in\cal N}\!\!N \hbox{ and }
\^G=\prod_{N\in\cal N}\!\!N.
$$

\Proof
$$
w(\underbrace{\^N,\^N,\dots,\^N}_{m-1\rm\;times},\^G,\^G,\dots,\^G)=
w(\underbrace{\^N,\^N,\dots,\^N}_{m-1\rm\;times},
\prod_{N\in\cal N}\!\!N,\^G,\dots,\^G)=
\prod_{N\in\cal N}\!\!w(\underbrace{\^N,\^N,\dots,\^N}_{m-1\rm\;times},
N,\^G,\dots,\^G).
$$
But $\^N\subseteq N$ and $\^G\subseteq G$, therefore, each factor of the
last product is contained in the group
$$
w(\underbrace{N,N,\dots,N}_{m\rm\;times},G,G,\dots,G),
\quad\hbox{
which is trivial by assumption.}
$$

\medskip
%\noindent
As a corollary, we obtain a strengthened version of the Khukhro--Makarenko
theorem with an explicit bound for the index.

\Th.
If a group $G$ contains a finite-index subgroup $N$
satisfying an outer commutator identity $w(x_1,\dots,x_t)=1$, then
$G$ contains a finite-index
subgroup $H$ satisfying the same
identity which is characteristic and even invariant under all
surjective endomorphisms.
In addition,
$$
\log_2|G{\,:\,}H|\le f^{t-1}(\log_2{|G{\,:\,}N|}) 
\quad\hbox{if the subgroup $N$ is normal,}
\eqno{(1)}
$$
and, therefore,
$\log_2|G{\,:\,}H|\le f^{t-1}(\log_2{|G{\,:\,}N|!})$ in the general case,
where $f^k(x)$ is the $k$-th iteration of the function $f(x)=x(x+1)$.

\Proof
For simplicity, we give a construction of a characteristic subgroup.
To obtain a subgroup invariant under all surjective
endomorphisms, one should replace all automorphisms by surjective
endomorphisms in the argument below.

Consider the subgroup
$
G_1=\!\prod\limits_{\varphi\in{\Aut G}}\!\!\!\!\!\varphi(N).
$
This subgroup is characteristic, and
$|G{\,:\,}G_1|\le |G{\,:\,}N|$. Clearly, $G_1$
is a product of at most $\log_2 |G{\,:\,}N|+1$ automorphic images of $N$
(because the chain
$N\subseteq N\phi_1(N)\subseteq N\phi_1(N)\phi_2(N)\subseteq\dots$
cannot have more than $\log_2 |G{\,:\,}N|+1$ different subgroups). Thus,
$$
G_1=\prod_{k=0}^{p_1}{\varphi'_k(N)},
\quad\hbox{where $\varphi'_k \in \Aut G$ and $p_1\le l_0\:=\log_2 |G{\,:\,}N|$}.
$$
Now, consider the subgroup
$
N_1=\bigcap\limits_{k=0}^{p_1}{\varphi'_k(N)}.
$
The index of an
intersection of subgroups does not exceed the product of their
indices (see, e.g., [KaM82]); hence,
$$
l_1\:=\log_2 |G{\,:\,}N_1|\le \log_2(|G{\,:\,}N|^{p_1+1})= (p_1+1)l_0\le (l_0+1)l_0=
f(l_0).
$$
By Lemma, we have
$$
w(N_1,\dots,N_1,G_1)=1.
$$

The next step is to consider the subgroups
$$
G_2=\!\prod_{\varphi\in{\Aut G_1}}\!\!\!\!\!{\varphi(N_1)}=
\prod_{k=0}^{p_2}{\varphi''_k(N_1)}
\quad\hbox{and}\quad
N_2=\bigcap_{k=0}^{p_2}{\varphi''_k(N_1)},
\quad\hbox{where }
\varphi''_k \in \Aut G_1
\hbox{ and }
p_2 \le \log_2 |G{\,:\,}N_1|= l_1 \le f(l_0).
$$
Clearly, $G_2$ is characteristic in $G$ (and even in $G_1$),
$$
\log_2 |G{\,:\,}G_2|\le \log_2 |G{\,:\,}N_1|=l_1 \le f(l_0),
\hbox{ and }
l_2\:=\log_2|G{\,:\,}N_2|\le \log_2(|G{\,:\,}N_1|^{p_2+1}) = 
(p_2+1)l_1\le f(l_1) \le f(f(l_0)).
$$
By Lemma, we have
$$
w(N_2,\dots,N_2,G_2,G_2)=1.
$$

Continuing in the same manner, at the $t$-th step, we obtain
a characteristic subgroup
$$
G_t=\!\prod_{\varphi\in{\Aut G_{t-1}}}\!\!\!\!\!{\varphi(N_{t-1})}=
\prod_{k=0}^{p_t}{\varphi_k^{(t)}(N_{t-1})},
\quad\hbox{where }
\varphi_k^{(t)}\in\Aut G_{t-1},
$$
such that
$$
w(G_t,\dots,G_t)=1
\quad\hbox{and}\quad
\log_2|G{\,:\,}G_t| \le \log_2|G{\,:\,}N_{t-1}|=l_{t-1} \le f(l_{t-2}) \le
f(f(l_{t-3}))\le 
\dots\le f^{t-1}(l_0).
$$
The subgroup $H=G_t$ is as required, which proves the theorem.

%\Замечание 1.
%Наше доказательство остаётся верным и в случае, когда
%индекс подгруппы $N$ является бесконечным кардиналом. В этой ситуации
%оценку (1) следует понимать так: $|G{\,:\,}H|\le|G{\,:\,}N|$. Таким образом,
%например, всякая группа, содержащая нормальную центрально метабелеву
%подгруппу счётного индекса, содержит также характеристическую центрально
%метабелеву подгруппу индекса не больше континуума.

\Remark.
In papers [KhM07] and [MKh07], there is no explicit bound,
but it is mentioned that an estimate can be obtained from the proof.
According to our calculations, inequality (1) is better.

\REFERENCES

\[KaM82]
Kargapolov M.I., Merzlyakov Yu.I.
Fundamentals of group theory.
Moscow: ``Nauka", 1982.

\[MKh07]
Makarenko N.Yu., Khukhro E.I.
Large characteristic subgroups satisfying multilinear commutator identities
{// Dokl. Akad. Nauk.} 2007. V.412. no. 5. P.594--596.

\[KhM07]
Khukhro E.I., Makarenko N.Yu.
Large characteristic subgroups satisfying multilinear commutator identities
{// J. London Math. Soc.} 2007. {V.75}. no.3, P.635--646.

\end